\documentclass[10pt,epsfig]{amsart}
\usepackage{graphicx}
\newtheorem{thm}{Theorem}[section]

\newtheorem{lem}{Lemma}[section]

\theoremstyle{definition}

\theoremstyle{remark}
\newtheorem{rk}{Remark}[section]
\numberwithin{equation}{section}

\newcommand{\vertiii}[1]{{\left\vert\kern-0.25ex\left\vert\kern-0.25ex\left\vert #1 \right\vert\kern-0.25ex\right\vert\kern-0.25ex\right\vert}}

\def \Proof{\noindent{\it Proof.}\quad}

\def \x{{\bf x }}

\def \M{{\bf M }}

\def \k{{\bf k}}
\def \eb{{\bf e}}

\def \a{\alpha}
\def \O{\Omega}
\def \o{\omega}
\def \vph{\varphi}
\def \th{\theta}

\def \n{\nu}

\def \vphh{\hat{\vph}}
\def \Hm{\rm H}
\def \Lm{\rm L}

\def \nab{\nabla}
\def \div{{\rm div}}
\def \pd{\partial}
\def \({\big(}
\def \){\big)}
\def \<{\langle}
\def \>{\rangle}
\def \neq{\not=}
\def \rha{\rightarrow}

\def \faa{\forall}
\def \S{\widetilde{S}}
\def \C{\widetilde{C}}

\def \Rbb{\mathbb{R}}

\begin{document}
\title[]{Numerical dispersion analysis of the convected Helmholtz equation}

\author{Ohsung Kwon and Imbo Sim}
\thanks{National Institute for Mathematical Sciences,70, Yuseong-daero 1689 beon-gil, Yuseong-gu, Daejeon, Republic of Korea ({imbosim@nims.re.kr}).}

 \dedicatory{}

\begin{abstract}
We present the numerical dispersion effects in solving the convected Helmholtz equation by the conforming and nonconforming quadrilateral finite elements. Particularly, we evaluate the dispersion relations for the numerical schemes.
The dispersive behaviors are analyzed by focusing on the Mach number and the angular frequency. Numerical experiments are conducted to verify the relations between the numerical dispersions and the computational errors.
\end{abstract}
\maketitle

\date{}

\section{Introduction}

In this paper, we consider the convected Helmholtz equation in the presence of a mean flow.
This equation is generated from the linearized Euler equations by reducing it for the pressure field  to describe a propagating wave to the mean flow (see \cite{PS}). 
So it has been applied in sciences and engineering problems, for instance, in aeroacoustics (\cite{Go, Ho}). In addition, various studies have been devoted as follows:
B{\'e}cache $\it et\, al.$ used perfectly marched layers for the convected Helmholtz equation to design efficient numerical absorbing boundary conditions in \cite{BBL}.
Casenave $\it et\, al.$ in \cite{CES} computed a linear acoustic wave propagation at a fixed frequency in the presence of flow using the coupled BEM-FEM.
An algebraic subgrid scale finite element method was presented by Guasch $\it et\, al.$ to improve the accuracy of the Galerkin finite element solution in \cite{GC}.
However, issues concerning the dispersion properties of the finite element method remain unsolved yet.

The dispersion relation concerns the angular frequency of a wave to the wavenumber and the Mach number. From this relation the phase and group velocities of the wave are derived.  Sometimes they lead instabilities due to the opposite signs (see \cite{Ar, BFJ, RH}).
Hence the dispersion analysis is an important issue. For the Helmholtz equation, there are numerous studies concerning the dispersion properties. Harari $\it et\, al.$ in \cite{HH1,HH2} used the Galerkin least squares to solve the Helmholtz equation and \cite{IB1,IB2,IB3} developed the generalized finite element method. Other approaches are appending the element boundary residuals to the Galerkin approximation in \cite{OP} and the nonconforming element in \cite{ZGS}. 

The aim of this paper is to analyze the numerical dispersion behaviors of the convected Helmholtz equation using the conforming and nonconforming finite element methods. 
Particularly, we investigate the difference of the continuous and numerical angular frequencies $\o^h-\o$ and fine that it relates to the numerical errors.


\section{Numerical dispersion relation of the convected Helmholtz equation}

Let us consider the convected Helmholtz equation in a uniform mean flow on the domain $\O\in \Rbb^2$:
\begin{equation}\label{ConvHelm2}
\begin{aligned}
-\div(A\nab p)-2i\o\, \M\cdot\nab p-\o^2 p&=0\mbox{ in } \O,\\
\n\cdot(A\nab p)-i\o p&=g\mbox{ on }\pd\O.
\end{aligned}
\end{equation}
Here, $\M=(M,0)$ for the Mach number $M$ and $A$ is a $2\times2$ diagonal matrix with $A_{11}=1-M^2$ and $A_{22}=1$.
%
Then the variational form of the problem \eqref{ConvHelm2} is to find $p\in\Hm^1(\O)$ such that
\begin{equation}\label{WK}
\(A\nab p,\nab v\)-2i\o (\M\cdot\nab p, v\)-\o^2\(p, v\)-i\o\<p,v\>=\<g, v\>, \quad\faa v\in\Hm^{1}(\O),
\end{equation}
where $(\cdot,\cdot)$ and $\<\cdot,\cdot\>$ are the $\Lm^2$-inner products such that $(u,v)=\int_{\O}uv\,dxdy$ and $\<u,v\>=\int_{\pd\O}uv\,ds$, respectively.

Next we will investigate the dispersion relations by the conforming and nonconforming finite element methods using quadrilateral elements of the lowest order. In particular, two numerical schemes are used: P1 conforming (P1-C) method, Rannacher-Turek nonconforming (RT-NC) method. Here, there are two types of RT elements. We set RT1 element by Rannacher-Turek element with the midpoints of edges as the degrees of freedom and RT2 element by one with mean integrals over edges (see \cite{RT}).

\begin{figure}[h!]
\begin{center}
\includegraphics[height=4.5cm]{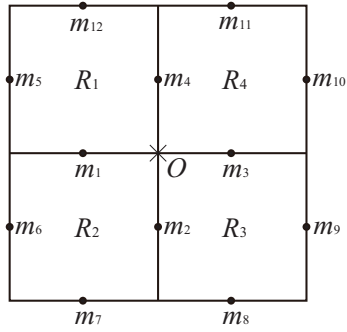}
\caption{Computational region $\O=[-h,h]^2$}
\label{FIG1}
\end{center}
\end{figure}

We first evaluate the dispersion relation by RT1-NC method. To do so, let $p^h$ be the numerical solution of the problem \eqref{WK} such that
\begin{equation}\label{ph}
p^h=\exp\{i(k_1^h x+k_2^hy)\},
\end{equation}
which is a plane wave propagating with the numerical wave vector $\k^h=(k_1^h,k_2^h)$.
Since the plane waves have the same structure on each quadrilateral mesh, we restrict the domain $\O=[-h,h]^2$ as Figure \ref{FIG1}. Then $\O$ contains 12 midpoints $m_j$ and 4 rectangles $R_l$.
The central point $O$ is not used in RT1-NC method, but it needs to define a global test function $\vph_G$. In particular, $\vph_G$ is defined by
\begin{equation*}\label{GB}
\vph_G=\left\{
\begin{aligned}
&\vphh_b|_{R_1}+\vphh_t|_{R_2}\mbox{ on } m_1,\\
&\vphh_r|_{R_2}+\vphh_l|_{R_3}\mbox{ on } m_2,\\
&\vphh_b|_{R_4}+\vphh_t|_{R_3}\mbox{ on } m_3,\\
&\vphh_r|_{R_1}+\vphh_l|_{R_4}\mbox{ on } m_4,
\end{aligned}
\right.
\end{equation*}
where $\vphh_j (j=r,l,t,b)$ is a local basis function on $R_l$ and the subscript $r$, $l$, $t$, and $b$ mean the right, left, top, and bottom on $R_l$, respectively.
Similarly, we denote $\vph_j$ the basis function at $m_j$ for $j=1,...,12$, for examples, $\vph_1=\vphh_b|_{R_1}+\vphh_t|_{R_2}$.
Using these basis functions $\vph_j$ of RT1-NC method, the numerical solution $p^h$ is represented by $p^h=\sum_{j=1}^{12}p_j\vph_j$ for constants $p_j$.
Inserting $p^h$ and $\vph_G$ into \eqref{WK} and using the middle point rule to approximate the inner products on the boundary, we have
\begin{equation*}
\sum_{j=1}^{12}p_j\big\{(A\nab\vph_j,\nab\vph_G)-2i\o (\M\cdot\nab\vph_j, \vph_{G})-\o^2(\vph_j, \vph_G)\big\}=0.
\end{equation*}
By the direct calculation of the inner products, we have
\begin{equation*}\label{PPP}
\begin{aligned}
&-\frac{\o^2h^2}{24}\(10\sum_{j=1}^4p_j+\sum_{j=5}^{12}p_j\)\\
&-\frac{iM\o h}{3}\big\{2(p_3-p_1)+2(p_{10}-p_6)+2(p_9-p_5)+(p_8-p_{12})+(p_{11}-p_{7})\big\}\\
&+2\sum_{j=1}^4p_j-\sum_{j=5}^{12}p_j+M^2\big\{p_5+p_6+p_9+p_{10}-2(p_2+p_4)\big\}=0.
\end{aligned}
\end{equation*}
By \eqref{ph}, $\sum_{j=1}^{12}p_j\vph_j=\exp\{i(k_1^h x+k_2^hy)\}$. Using the properties $\vph_j(m_j)=1$ and $\vph_j(m_l)=0$ for $j\neq l$, $p_j$ is expressed in terms of $k_j^h$ and $h$, e.g., $p_1=\exp(-i k_1^hh/2)$.
So, letting $C_j=\cos(k_j^hh/2)$ and $S_j=\sin(k_j^hh/2)$ for $j=1,2$, the dispersion relation of RT1-NC method is given by

\begin{equation}\label{O_RT1}
\begin{aligned}
\o=\frac{4M}{h}\Bigg(G_1+\sqrt{G_1^{\,2}
 +\frac{3}{2}\frac{(C_1+C_2)(1-C_1C_2)-M^2S_1^2C_2}{M^2(C_1+C_2)(2+C_1C_2)}}\Bigg),
\end{aligned}
\end{equation}
where
$$
G_1=\frac{S_1C_2(2C_1+C_2)}{(C_1+C_2)(C_1C_2+2)}.
$$
Let $\C_j=\cos(k_j^hh)$, $\S_j=\sin(k_j^hh)$ for $j=1,2$.
By the same method used in \eqref{O_RT1}, we obtain the dispersion relations for RT2-NC and P1-C methods such that
\begin{equation}\label{O_RT2}
\begin{aligned}
\o=&\frac{6M}{h}\Bigg(G_2+\sqrt{G_2^{\,2}+\frac{2}{3}\frac{k_1^h\S_2(1-\C_1)+k_2^h\S_1(1-\C_2)-M^2k_1^h\S_2(1-\C_1)}{M^2\(k_1^h\S_2(5+\C_1)+k_2^h\S_1(5+\C_2)\)}}\Bigg),\\
\o=&\frac{3M}{h}\Bigg(G_3+\sqrt{G_3^{\,2}-\frac{2}{3}\frac{2\C_1\C_2+\C_1+\C_2-4+M^2(2-\C_1\C_2-2\C_1+\C_2)}{M^2\(\C_1\C_2+2(\C_1+\C_2)+4\)}}\Bigg), 
\end{aligned}
\end{equation}
where
\begin{equation*}
G_2=\frac{k_1^h\S_1\S_2+k_2^h(1-\C_1)(1+\C_2)}{k_1^h\S_2(5+\C_1)+k_2^h\S_1(5+\C_2)},\quad G_3=\frac{\S_1(2+\C_2)}{\C_1\C_2+2(\C_1+\C_2)+4},
\end{equation*}
respectively.

\begin{rk}
The dispersion relation of the continuous problem \eqref{ConvHelm2} is of the form:
\begin{equation}\label{DR1}
\o=k_1M+|\k|,
\end{equation}
which is derived from the characteristic equation of \eqref{ConvHelm2} by \eqref{ph}.
\end{rk}

\section{Analysis of the numerical dispersion relation}

\begin{figure}
\begin{center}
\includegraphics[height=12cm]{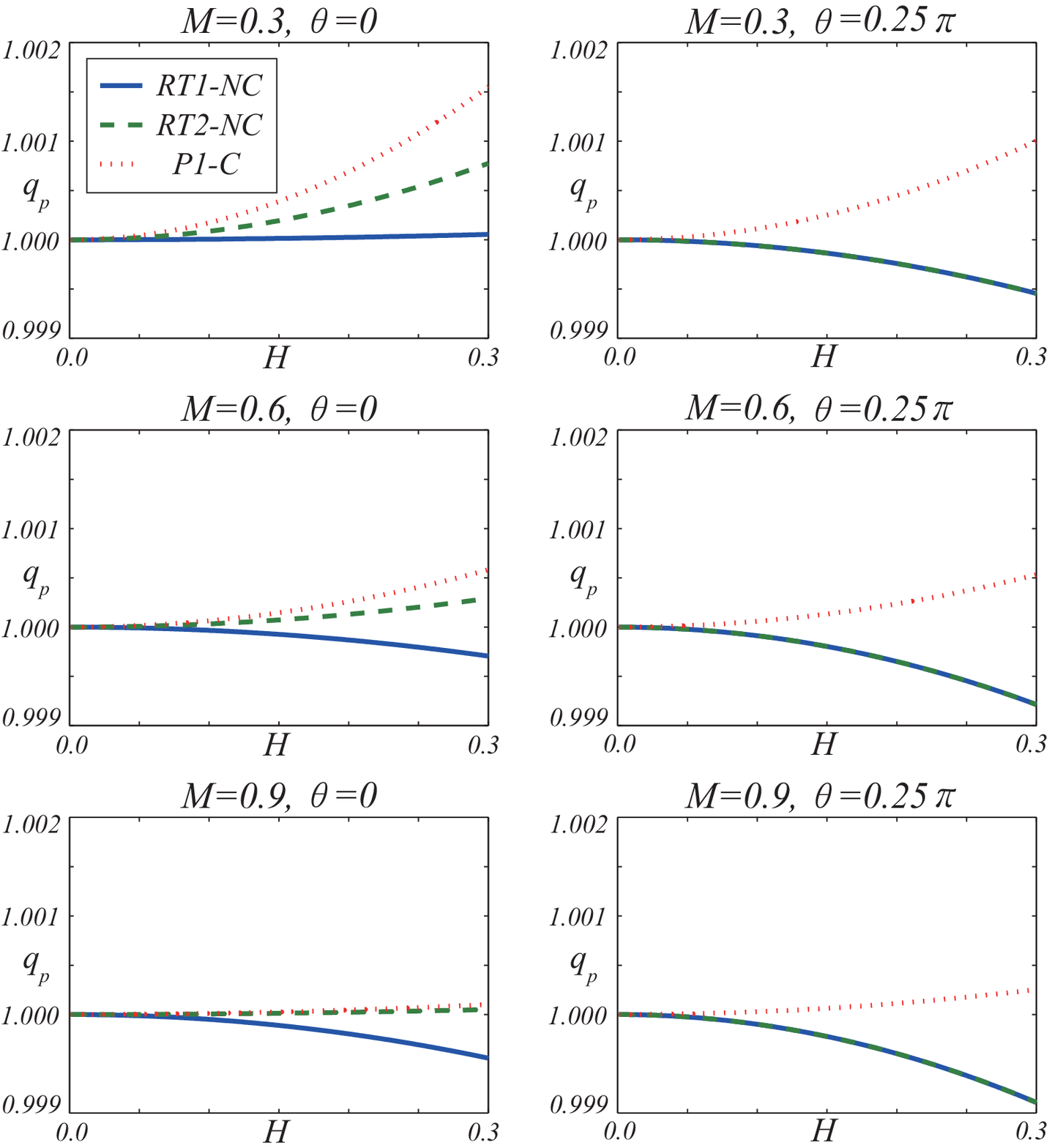}
\caption{Dispersion quotient $q_p$ bias $H$ for $M=0.3, 0.6, 0.9$ and $\th=0, 0.25\pi$}
\label{FIG2-1}
\end{center}
\end{figure}
\begin{figure}
\begin{center}
\includegraphics[height=12cm]{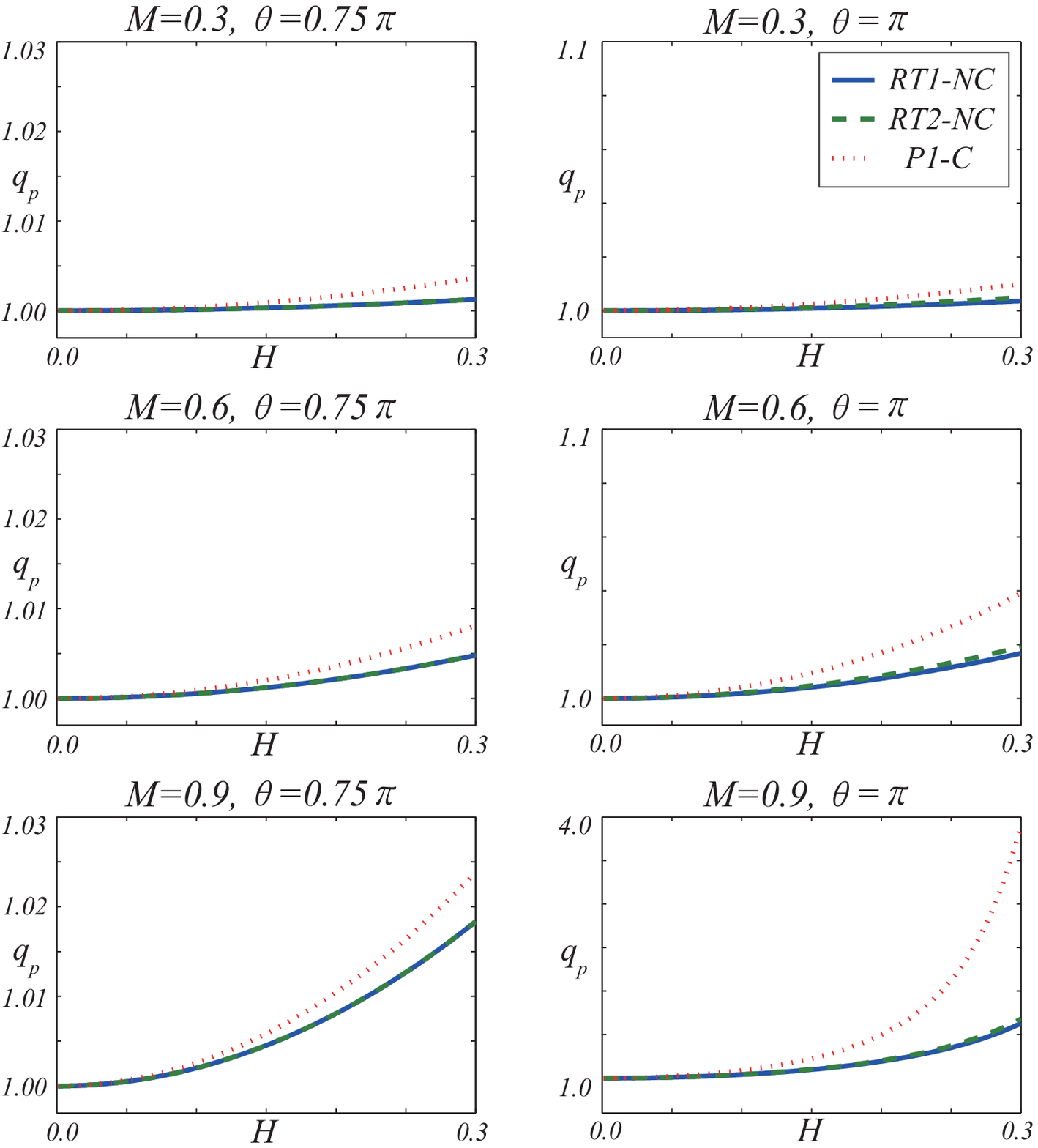}
\caption{Dispersion quotient $q_p$ bias $H$ for $M=0.3, 0.6, 0.9$ and $\th=0.75\pi, \pi$}
\label{FIG2-2}
\end{center}
\end{figure}

To derive the dispersive behaviors of \eqref{O_RT1}--\eqref{O_RT2}, we define a new variable $\th$ that is the angle for the direction of the wave propagation so that $\k=k(\cos\th,\sin\th)$ with $k=|\k|$. Then \eqref{DR1} is represented by $\o=k(1+M\cos\th)$. As a result, the dispersion relations in \eqref{O_RT1}--\eqref{O_RT2} are reformulated in terms of $\o^h$. Then we can expand $\o$ in terms of $h$ such that
\begin{equation}\label{TE}
\o= \o^h\Big(1+\sum_{j=1}^\infty A_j(M,\th)\,(\o^hh)^{2j}\Big),
\end{equation}
where $A_j(M,\th)$ is a function depending on $M$ and $\th$.
This yields that $\o^h$ converges to the continuous angular frequency $\o$ as $h\rha 0$.

To compare the dispersive behaviors of the numerical schemes, we define two dispersion quotients in \cite{Co}:
$$
q_p=v_p/v_p^h,\quad q_g=v_g/v_g^h,
$$
where $v_p$ and $v_g$ are the phase and group velocities in the direction of the numerical wave vector $\k^h$ such that $v_p=\o/k^h$ and $v_g=|\pd\o/\pd\k^h|$, respectively. Similarly, the numerical phase and group velocities are $v_p^h=\o^h/k^h$ and $v_g^h=|\pd\o^h/\pd\k^h|$, respectively. Then $q_p$ and $q_g$ measure the errors in the phase and group velocities.

\begin{figure}
\begin{center}
\includegraphics[height=12cm]{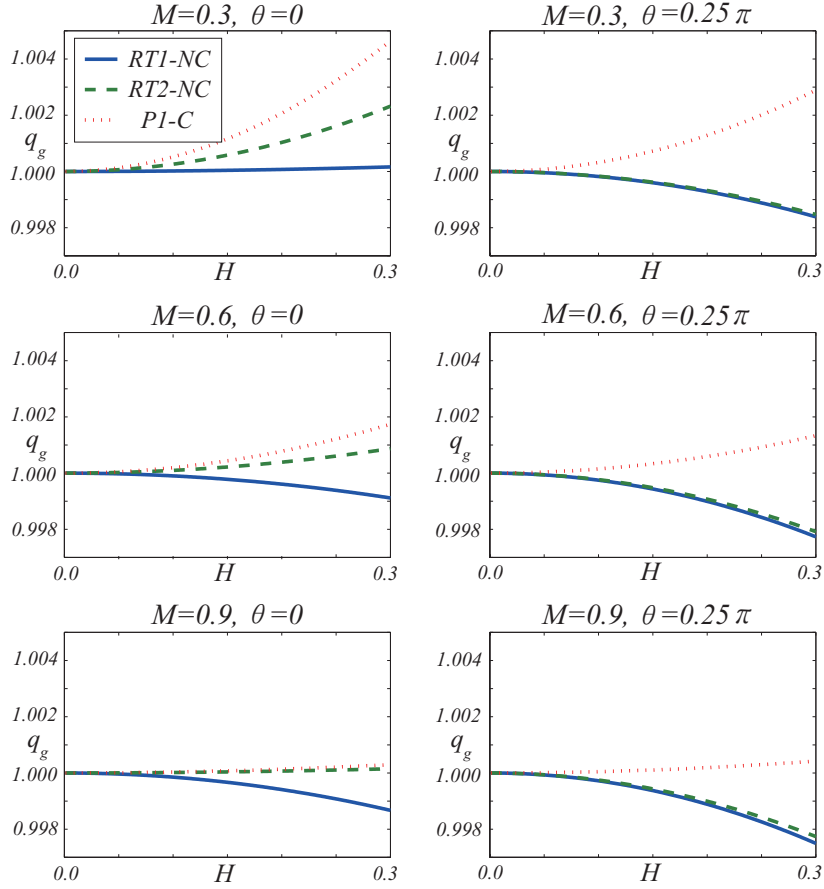}
\caption{Dispersion quotient $q_g$ bias $H$ for $M=0.3,0.6,0.9$ and $\th=0, 0.25\pi$}
\label{FIG3-1}
\end{center}
\end{figure}

\begin{figure}
\begin{center}
\includegraphics[height=12cm]{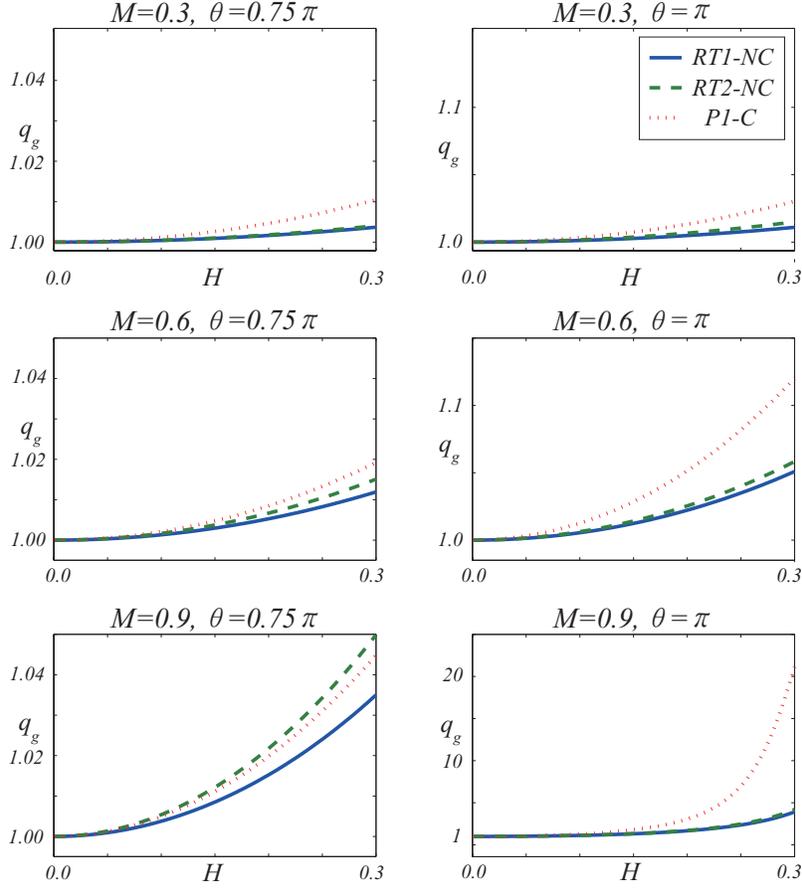}
\caption{Dispersion quotient $q_g$ bias $H$ for $M=0.3,0.6,0.9$ and $\th=0.75\pi, \pi$}
\label{FIG3-2}
\end{center}
\end{figure}

Using \eqref{TE}, we have $q_p=1+\sum_{j=1}^\infty A_j(M,\th)\,(\o^hh)^{2j}$. To observe the effects of $M$ and $\th$ clearly, we set $H:=\o^hh\in[0,0.3]$, since the numerical experiments will be done on this interval.
Since the mean flow passes horizontally, the dispersion quotient $q_p$ is symmetric with respect to $x$ axis, so it suffices to consider $\th\in[0,\pi]$. In particular, we use four angles of $\th$ such as $0, 0.25\pi, 0.75\pi,$ and $\pi$,
since the different behaviors of the numerical schemes are well observed on these angles. For Mach number $M$, we set $M=0.3, 0.6,$ and $0.9$. The asymptotic behavior of $q_g$ is also similar to $q_p$, so we use the same setting.
Under these conditions, the dispersion quotients $q_p$ and $q_g$ are illustrated as Figures \ref{FIG2-1}--\ref{FIG3-2}, respectively.


Figures \ref{FIG2-1}-\ref{FIG3-2} show the effects of $M$ and $\th$ to the dispersion quotients of numerical schemes. For RT-NC method, the dispersion quotients get away to 1 as $M$ and $\th$ are larger. Specially, it is significantly large at $M=0.9$ and $\th=\pi$. For P1-C method, it shows the anisotropy behavior to $M$. For large $\th$, it appears the similar behaviors to RT-NC method. Meanwhile, the dispersion quotients approach to 1 as $M$ is larger for small $\th$.
This means that P1-C method is more efficient than RT-NC method for small $\th$ and large $M$, though it does not for other cases.

In fact, the behaviors of the dispersion quotients are determined by the difference of the numerical and continuous angular frequencies. 
In next theorem, we investigate the behavior of the difference $\o^h-\o$ using techniques in \cite{BIPS,OP}.

\begin{thm} \label{THM1}
Let $\o$ and $\o^h$ be continuous and numerical angular frequencies, respectively. We define the dispersion error by $|\o^h-\o|$. Then $|\o^h-\o|=A_1(M,\th)\o^3h^2+O(\o^5 h^4)$, where
\begin{equation*}\label{Err}
A_1=
\left\{
\begin{aligned}
&\frac{|2(1+\cos(4\th))+4M(\cos(3\th)-3\cos\th)+M^2(\cos(4\th)-6\cos(2\th)-7)|}{384(1+M\cos\th)^3},\\
&\frac{|2(1+\cos(4\th))+4M(\cos(3\th)-\cos\th)+M^2(\cos(4\th)-2\cos(2\th)-3)|}{192(1+M\cos\th)^3},\\
&\frac{|3+\cos(4\th)-4M^2\cos^4\th|}{96(1+M\cos\th)^3},
\end{aligned}
\right.
\end{equation*}
which are the leading coefficient functions for RT1-NC, RT2-NC, and P1-C methods, respectively.
\end{thm}
\begin{figure}
\begin{center}
\includegraphics[height=8cm]{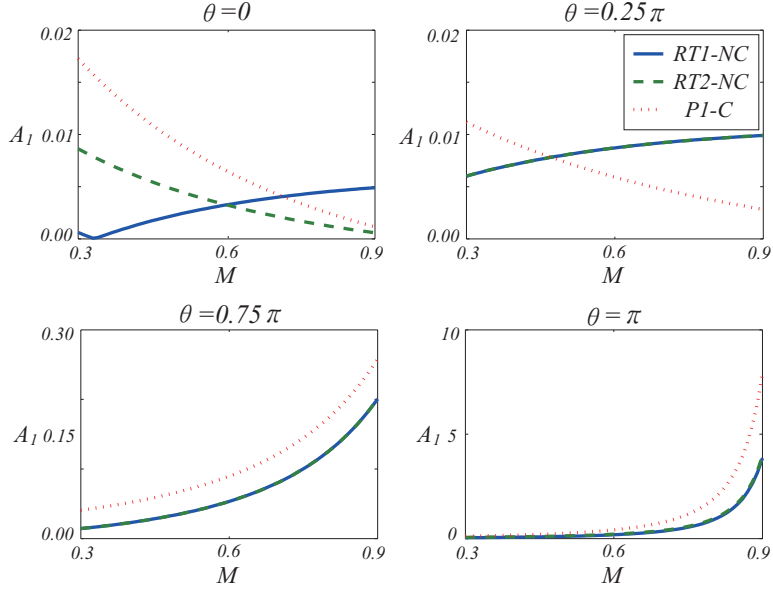}
\caption{Graphes of $A_1(M,\th)$ bias $M$ for $\th=0, 0.25\pi, 0.75\pi, \pi$}
\label{FIG4}
\end{center}
\end{figure}

Generally, if $h$ is sufficiently small so that $\o^3h^2\leq 1$ for large $\o$, then the main behavior of the dispersion error is determined by $A_1(M,\th)$.
This fact is confirmed by that the effects of $M$ and $\th$ to the dispersion quotients in Figures \ref{FIG2-1}--\ref{FIG3-2} follows the behaviors of $A_1(M,\th)$ in Figure \ref{FIG4}.

To check the effects of the dispersion errors to the computational errors of the numerical schemes, we define the norm:
\begin{equation}\label{norm}
\vertiii{p}^2:=\|\sqrt{A}\nab p\|_{L^2(\O)}^2+\o^2\|p\|_{L^2(\O)}^2,
\end{equation}
which is the energy norm of \eqref{WK} for the real part. To calculate the numerical solution, we need the boundary date $g$. It could be verified by $p^h=\exp\{i(k_1^h x+k_2^hy)\}$ such that
$$
g=i(k\n\cdot(A\eb_r)-\o)\exp(ik\x\cdot\eb_r),
$$
where $\eb_r=(\cos\th,\sin\th)$ and $\x=(x,y)$.

\begin{figure}
\begin{center}
\includegraphics[height=12cm]{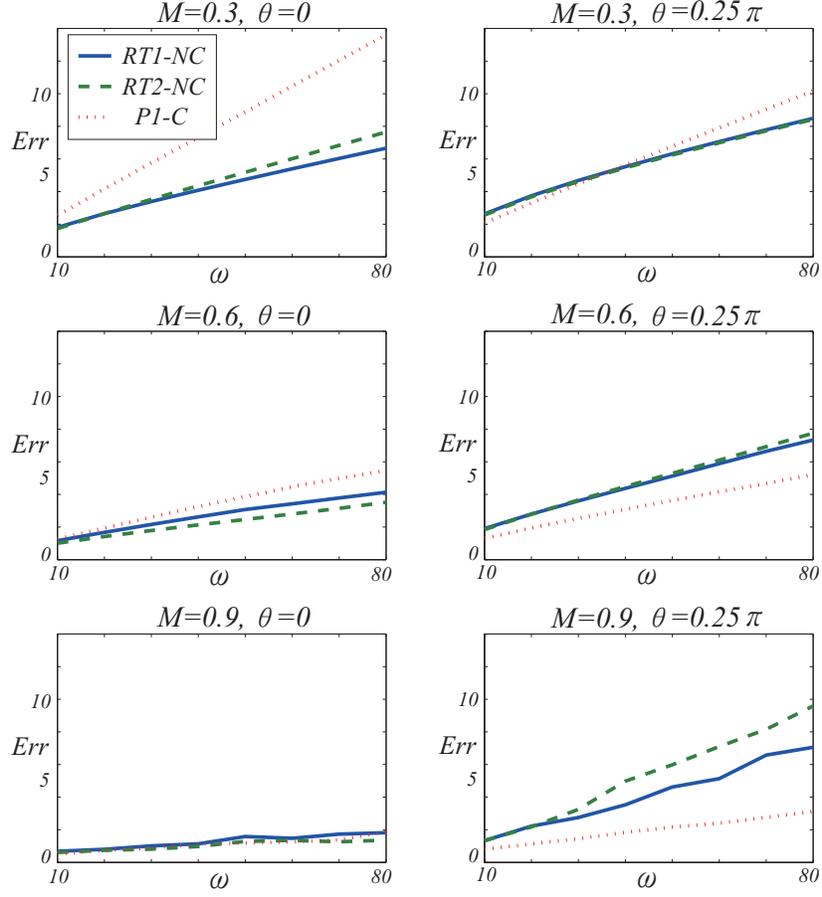}
\caption{Numerical errors bias $\o$ for $M=0.3,0.6,0.9$ and $\th=0,0.25\pi$}
\label{FIG5-1}
\end{center}
\end{figure}

\begin{figure}
\begin{center}
\includegraphics[height=12cm]{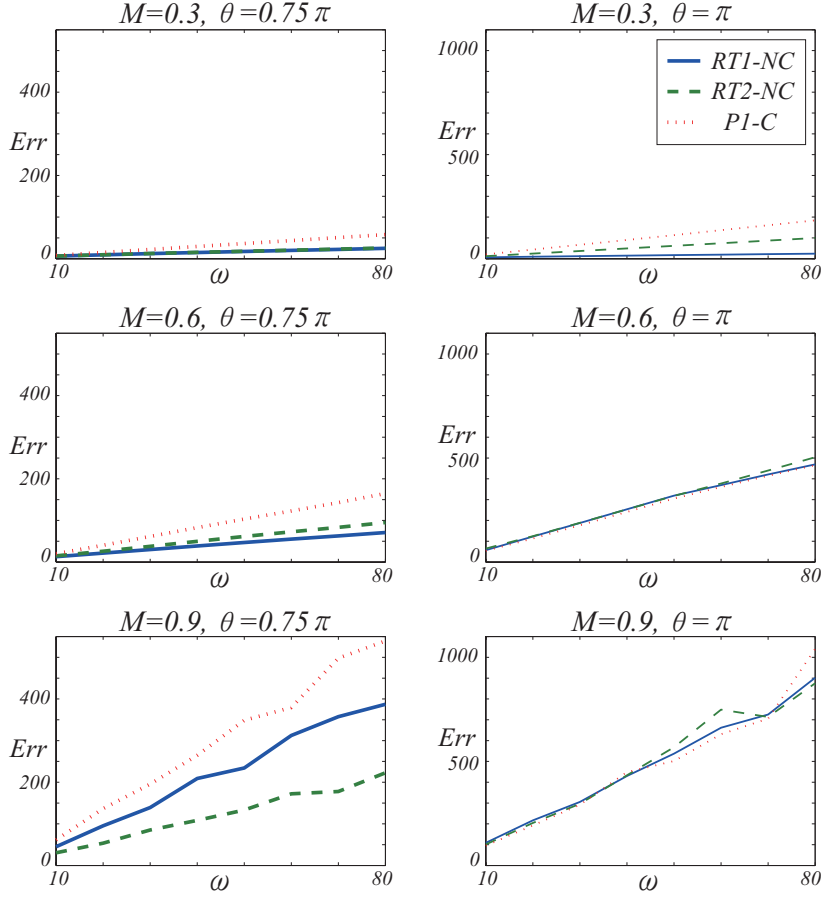}
\caption{Numerical errors bias $\o$ for $M=0.3,0.6,0.9$ and $\th=0.75\pi,\pi$}
\label{FIG5-2}
\end{center}
\end{figure}

In order to observe the effect of $A_1(M,\th)$ clearly, we set $\o^3h^2=1$ and solve the problem by growing $\o$ from $10$ to $80$. 
Figures \ref{FIG5-1} and \ref{FIG5-2} illustrate the numerical error $Err:=\vertiii{p-p^h}$ for $M$ and $\th$, which yields that the error behaviors nearly follow $A_1(M,\th)$ in Figure \ref{FIG4}.


\section{Comparison to the Helmholtz formulation} 

The convected Helmholtz equation in \eqref{ConvHelm2} could be reformulated as the Helmholtz equation by the following lemma.
\begin{lem}\label{LEM4.1}
Let $\a(x)=e^{i\o M x/(1-M^2)}$. If we set $u(x,y):=p(x,y)\a(x)$, then $u(x,y)$ is a solution of the problem:
\begin{equation}\label{Helm}
-\div(A\nab u)-\frac{\o^2}{1-M^2}u=0\quad\mbox{ in }\O.
\end{equation}
\end{lem}
\Proof
Let $d=1-M^2$. Inserting $u(x,y)=p(x,y)\a(x)$ into \eqref{Helm}, we have
$$
\begin{aligned}
-\div(A\nab u)-\frac{\o^2}{1-M^2}u&=-d(p\a)_{xx}-(p\a)_{yy}-\o^2p\a/d\\
&=-d(p_{xx}\a+2p_x\a_x+p\a_{xx})-p_{yy}\a-\o^2p\a/d\\
&=-d\a\(p_{xx}+2i\o Mp_x/d-M^2\o^2p/d^2\)-p_{yy}\a-\o^2p\a/d\\
&=\a\(-d p_{xx}-p_{yy}-2i\o Mp_x-\o^2p\).
\end{aligned}
$$
\qed

Using Lemma \ref{LEM4.1}, we can find the solution $p$ of the problem \eqref{ConvHelm2} without considering the convection term by solving \eqref{Helm}.
However the problem \eqref{Helm} has a stability deterioration for large Mach number by the term $1/(1-M^2)$ in \eqref{Helm}. 
This phenomena is generated by the dispersion error $\o^h-\o$ as follows:

We first calculate the dispersion relations of the problem \eqref{Helm}.
By \eqref{ph} and \eqref{DR1}, the numerical solution $u^h$ of $u$ is of form $u^h=\exp\{i(k_3^hx+k_2^hy)\}$ with $k_3^h=(k^h_1+k^hM)/(1-M^2)$. Using the process stated in Section 2, the dispersion relations for RT1-NC, RT2-NC, and P1-C methods are given by
$$
\o=\left\{\begin{aligned}
&\frac{2}{h}\sqrt{6(1-M^2)\frac{(C_2+C_3)(1-C_2C_3)+M^2C_2(C_3^2-1)}{(C_2+C_3)+(2+C_2C_3)}},\\
&\frac{2}{h}\sqrt{6(1-M^2)\frac{k_2\S_3(1-\C_2)+k_3\S_2(1-\C_3)-M^2k_3\S_2(1-\C_3)}{k_2\S_3(\C_2+5)+k_3\S_2(\C_3+5)}},\\
&\frac{1}{h}\sqrt{6(1-M^2)\frac{4-2\C_2\C_3-\C_2-\C_3+M^2(\C_2\C_3-\C_2+2\C_3-2)}{\C_2\C_3+2(\C_2+\C_3)+4}},
\end{aligned}\right.
$$
where $C_3=\cos(k_3^hh/2)$ and $\C_3=\cos(k_3^hh)$.
Under the same assumption to $h$ and $\o$ in Section 3, the main behavior of the dispersion error is determined by the leading coefficient $A_1(M,\th)$, so it suffices to consider $A_1(M,\th)$ for each numerical scheme. 
For RT1-NC method, $A_1(M,\th)$ is given by
$$
\begin{aligned}
A_1(M,\th)=&\bigg|\frac{(4+10M^2+28M^4-7M^6)+4M(4+11M^2-M^4)\cos\th}{{768(1-M^2)^2(1+M\cos\th)^4}}\\
&+\frac{4M^2(11-6M^2+2M^4)\cos(2\th)+4M(4-3M^2+M^4)\cos(3\th)}{768(1-M^2)^2(1+M\cos\th)^4}\\
&+\frac{(4-6M^2+4M^4-M^6)\cos(4\th)}{768(1-M^2)^2(1+M\cos\th)^4}\bigg|.
\end{aligned}
$$
For RT2-NC method, it is the double of one for RT1-NC method. For P1-C method, it is
$$
A_1=\frac{|\cos^4\th+4M\cos^3\th+6M^2\cos^2\th+4M^3\cos\th+M^4+(1-M^2)^3\sin^4\th|}{24(1-M^2)^2(1+M\cos\th)^4}.
$$
All $A_1(M,\th)$ blow up as $M\rha 1$, so it leads to instabilities for large $M$.

\section{Conclusion}
We have analyzed the dispersive behaviors of the convected Helmholtz equation by the conforming and nonconforming finite element methods.
Particularly, the dispersion relations of the numerical schemes are derived on the quadrilateral mesh and it is shown that the numerical dispersions converge to the continuous ones as $h\rha0$.
We also observed the effects of the dispersion error $|\o^h-\o|$ to the dispersion quotients and the numerical errors.

Consequently, our result informs the effects of the Mach number $M$ and the angular frequency $\o$ to the numerical errors.
So it provides the guidelines for the selection of an appropriate mesh size in terms of $M$ and $\o$ in solving numerically.


\end{document}